\documentclass{article}
\usepackage{amsthm}
\usepackage[T1]{fontenc}
\usepackage[frenchb, english]{babel}
\usepackage[hmargin=4cm,vmargin=4cm]{geometry}
\usepackage{amssymb}
\usepackage{amsmath}
\usepackage{pstricks}

\def\aro{{A}^{\!\!\!\raise5pt\hbox{$\scriptstyle \circ$}}}
\def\aroi{{A}^{\!\!\!\raise4pt\hbox{$\scriptscriptstyle \circ$}}}
\def\cro{\smash{{C}^{\!\!\!\raise5pt\hbox{$\scriptstyle \circ$}}}}
\def\croi{\smash{{C}^{\!\!\!\raise4pt\hbox{$\scriptscriptstyle \circ$}}}}
\begin{document}

\newtheorem{Th}[subsubsection]{Théorème}
\newtheorem{Pro}[subsubsection]{Proposition}
\newtheorem{De}[subsubsection]{Définition}
\newtheorem{Prt}[subsubsection]{Propriété}
\newtheorem{Prts}[subsubsection]{Propriétés}
\newtheorem{Le}[subsubsection]{Lemme}
\newtheorem{Hyp}[subsubsection]{Hypothèse}
\newtheorem{Cor}[subsubsection]{Corollaire}

%\newenvironment{desc}[1]{
%	\begin{list}{}
%		{\renewcommand{\makelabel}[1]{
%			\textup{\textsf{##1}}\hfill}
%		\settowidth{\labelwidth}{\textsf{#1 }}
%		\leftmargin=\labelwidth
%		\advance \leftmargin\labelsep
%}}
%	{\end{list}}

%\bigskip

\title{A proof of the equivalence of ensembles\\
for asymptotically decoupled fields\\
relying on Mosco's theorem}
\date{March 18, 2011}
\author{Pierre Petit\\ \\
Universit\'e Paris Sud}
%\affil{Universit\'e Paris Sud}
\maketitle

\selectlanguage{english}

\begin{abstract}
We give a proof of the identification $s = -p^*$ between the entropy and the opposite of the Fenchel-Legendre transform of the pressure for asymptotically decoupled fields relying on Mosco's theorem. The proof is longer than the one based on convex-tension in \cite{CeP11} but it explores the links between large deviations theory and Mosco's theorem.
\end{abstract}

\selectlanguage{french}

\begin{abstract}
Nous démontrons l'identification $s = -p^*$ entre l'entropie et l'opposée de la transformée de Fenchel-Legendre de la pression pour les champs asymptotiquement découplés à partir du théorème de Mosco. La preuve est certes plus longue que celle de \cite{CeP11} qui repose sur la convexe-tension mais explore les liens entre la théorie des grandes déviations et le théorème de Mosco.
\end{abstract}

\medskip

\textbf{Remarque :} La prochaine version de ce texte sera en langue anglaise.

\section{Introduction}
%\addcontentsline{toc}{section}{Introduction}

La démonstration exposée ici de l'égalité
$$
s = -p^*
$$
dans le cas asymptotiquement découplé et pour un espace de Fréchet séparable (plus précisément pour une limite projective d'espaces de Fréchet séparables dans lesquels il y a le découplage asymptotique) généralise, au cas asymptotiquement découplé, l'idée de \cite[chapitre 26]{Cer07} (reprise de \cite{Zab92a} et \cite{Zab92b}) d'utiliser le théorème de Mosco (\emph{cf.} \cite{Mos71}). Il existe une démonstration moins laborieuse reposant simplement sur la dualité convexe et l'utilisation de la convexe-tension (cf. \cite{CeP11}), mais la preuve présentée ici jette un nouveau pont entre grandes déviations et théorème de Mosco, ces liens étant encore mal compris. Nous nous servons de l'égalité vraie si la coordonnée $\sigma(0)$ du champ $\sigma$ en $0$ prend ses valeurs dans un ensemble compact ; puis, pour l'étendre à une mesure quelconque, nous approximons $\sigma$ par une suite de champs dont la coordonnée en $0$ prend ses valeurs dans un ensemble compact.

%\medskip

%Le cadre et les notations n'ont pas été repris pour être conformes à ceux du chapitre précédent. On pourra noter quelques variations de définitions. Par exemple, on considère d'emblée une mesure sur $X^{\mathbb{Z}^d}$, et non un système compatible de mesures. Nous refaisons la démonstration du lemme sous-additif dans le cas asymptotiquement découplé, cette fois pour une mesure, non plus invariante par toutes les translations sur le réseau $\mathbb{Z}^d$, mais seulement invariante par les translations d'un sous-réseau $(\ell \mathbb{Z})^d$ de $\mathbb{Z}^d$ (car les mesures que nous manipulerons pour l'approximation ne seront pas invariantes par translation).

\medskip

On montre le résultat dans $X$, limite projective d'espaces $Y_i$. Il est raisonnable d'imposer aux $Y_i$ d'être des espaces de Fréchet séparables. Voici en effet les hypothèses apparaissant au cours de la démonstration : pour tous $z$, la loi de $\eta(z) = f_i(\sigma(z))$ est tendue, les ouverts de $Y_i$ sont des réunions dénombrables de pavés ouverts mesurables (c'est le cas si $\tau(Y_i)$ est à base dénombrable d'ouverts), $Y_i$ est un e.v.l.c. séparable et métrisable (pour appliquer le théorème de Mosco), $Y_i$ est tonnelé (pour appliquer le théorème de Banach-Steinhaus et pour l'argument final faisant intervenir un tonneau).

%A corriger éventuellement :
%\begin{itemize}
%\item définir $\Lambda'(m) := \Lambda(m + g(m) + \ell)$.
%\item si $\lambda \in Y_i^*$, mesurabilité de $\lambda$ ?
%\item a-t-on $\mathcal{C} \varsubsetneq \mathcal{C}(X)$ ?
%\item mettre partout $X^*$ et pas $X'$ (cf. après le PGD).
%\item FKG n'est vérifiée que pour les événements croissants, d'habitude\ldots
%\end{itemize}

%\[
%\left.
%\begin{array}{l}
%\textrm{métrisable complet $\Rightarrow$ de Baire (TG IX, p. 55)}\\
%\textrm{e.v.l.c. de Baire $\Rightarrow$ tonnelé (EVT III, p. 25)}\\
%\textrm{e.v.l.c. à base dénombrable de voisinages de $0$ $\iff$ métrisable (EVT II, p. 25)}\\
%\textrm{(autrement dit $\mathcal{C}_0(X)$ dénombrable $\iff$ métrisable.)}
%\end{array}
%\right\} \Rightarrow \textrm{1., 5., 6. $\iff Y_i$ e.v.l.c.m. de Baire}
%\]

%Mesurabilité des applications linéaires ?

%Topologie localement convexe la plus fine : topologie ayant pour système de voisinage de $0$ les convexes symétriques absorbants. Toute application linéaire d'un tel espace dans un autre evlc est automatiquement continue.

%Regarder de plus près l'hypothèse Cerf 16.2.

\section{Lemme sous-additif}\label{scadre}

Soient $(Y, \tau)$ un espace de Fréchet séparable, $\mathcal{B}$ sa tribu borélienne et $d \in \mathbb{N}^*$.

\subsection{Découplage asymptotique convexe}

Pour tout sous-ensemble fini $S$ de $\mathbb{Z}^d$, on notera $\mathcal{C}(Y^S)$ l'ensemble des ouverts convexes de $Y^S$ relativement à la topologie produit. Les éléments de $\mathcal{C}(Y^S)$ sont mesurables relativement à la tribu produit. On entendra par \emph{boîte} un sous-ensemble cubique de $\mathbb{Z}^d$, autrement dit un ensemble de la forme
$$
\Lambda(z ; n) := z + [0, n[^d \cap \mathbb{Z}^d
$$
où $n \in \mathbb{N}^*$. On notera $\Lambda(n) = \Lambda(0 ; n)$. On notera $\textup{dist}$ la distance associée à la norme $|\cdot|_\infty$ sur $\mathbb{Z}^d$. Si $\xi \in Y^{\mathbb{Z}^d}$ et $z \in \mathbb{Z}^d$, on notera $\xi(z)$ la coordonnée suivant $z$ de $\xi$. Plus généralement, pour tout sous-ensemble fini $S$ de $\mathbb{Z}^d$, on notera
$$
\xi_S := (\xi(z))_{z \in S}
$$
La définition fondamentale suivante est reprise de \cite{Pfi01}. On se donne un champ $\eta$ à valeurs dans $\smash{\big(Y^{\mathbb{Z}^d}, \mathcal{B}^{\otimes \mathbb{Z}^d}\big)}$ et deux applications, $g$ et $c : \mathbb{N}^* \rightarrow [0, + \infty[$, telles que
$$
\frac{g(n)}{n} \to 0 \quad \textrm{et} \quad \frac{c(n)}{|\Lambda(n)|} \to 0
$$
On dit que $\eta$ est \emph{asymptotiquement découplé inférieurement pour les convexes} (\emph{c.a.d.i.}) de paramètre $(g, c)$ (ou encore $(g, c)$-c.a.d.i.) si, pour tous $z \in \mathbb{Z}^d$ et $n \in \mathbb{N}$, pour tout ensemble fini $S$ de $\mathbb{Z}^d$, pour tous ouverts convexes mesurables $C \in \mathcal{C}(Y^{\Lambda(z ; n)})$ et $D \in \mathcal{C}(Y^S)$,
$$
\textup{dist} \big(S , \Lambda(z ; n)\big) > g(n)
\Longrightarrow \mathbb{P} \big( \eta_{\Lambda(z ; n)} \in C, \, \eta_S \in D \big) \geqslant e^{-c(n)}  \mathbb{P} \big( \eta_{\Lambda(z ; n)} \in C \big) \mathbb{P} \big( \eta_S \in D \big)
$$
De plus, si $\mu \in \mathcal{M}_1^+\smash{\big(Y^{\mathbb{Z}^d}, \mathcal{B}^{\otimes \mathbb{Z}^d}\big)}$, on dit que $\mu$ est c.a.d.i. de paramètre $(g, c)$ si c'est la loi d'un champ $\eta$ c.a.d.i. de paramètre $(g, c)$. Si $f : Y \rightarrow Z$, on notera sans ambiguïté
$$
f(\eta) := \big(f(\eta(z))\big)_{z \in \mathbb{Z}^d}
$$

\begin{Pro}[Stabilité affine]\label{scadimos}
Soit $\eta$ un champ c.a.d.i. de paramètre $(g, c)$. Si $Z$ est un espace de Fréchet séparable et $f : Y \rightarrow Z$ une application affine, continue et mesurable, alors $f(\eta)$ est c.a.d.i. également, de paramètre $(g, c)$.
\end{Pro}

\textbf{Démonstration :} Il suffit de remarquer que, pour toute boîte $\Lambda$, l'application $\phi : \xi_\Lambda \rightarrow \big(f(\xi(z))\big)_{z \in \Lambda}$ est continue et mesurable.

\subsection{Contrôle local}

Introduisons aussi une hypothèse pour contrôler chaque site conditionnellement au reste de la configuration. On note $\mathcal{C}_0(Y)$ une base de voisinages convexes ouverts de $0$ dans $Y$. On se donne un champ $\eta$ à valeurs dans $\smash{\big(Y^{\mathbb{Z}^d}, \mathcal{B}^{\otimes \mathbb{Z}^d}\big)}$ et deux applications, $t : \mathcal{C}_0(Y) \rightarrow ]0, + \infty[$ et $\alpha : \mathcal{C}_0(Y) \rightarrow ]0, 1]$. On dit que $\eta$ est \emph{contrôlé localement} (\emph{c.l.}) de paramètre $(t, \alpha)$ (ou encore $(t, \alpha)$-c.l.) si, pour tout $V \in \mathcal{C}_0(Y)$, pour tout $z \in \mathbb{Z}^d$ et tout sous-ensemble fini $S$ de $\mathbb{Z}^d$, pour tout $D \in \mathcal{C}(Y^S)$,
$$
\mathbb{P} \big(\eta(z) \in t(V) \cdot V \, ; \, \eta_S \in D\big) \geqslant \alpha(V) \, \mathbb{P} (\eta_S \in D)
$$
De plus, si $\mu \in \mathcal{M}_1^+\smash{\big(Y^{\mathbb{Z}^d}, \mathcal{B}^{\otimes \mathbb{Z}^d}\big)}$, on dit que $\mu$ est c.l. de paramètre $(t, \alpha)$ si c'est la loi d'un champ $\eta$ c.l. de paramètre $(t, \alpha)$. Si $V$ est un voisinage convexe de $0$, sa jauge (ou fonctionnelle de Minkowski) est l'application $M_V : Y \rightarrow [0, +\infty]$ définie par
$$
\forall y \in Y \qquad M_V(y) = \inf \{ t \geqslant 0 ; y \in tV \}
$$

\begin{Pro}[Stabilité affine]\label{sclmos}
Soit $\eta$ un champ c.l. de paramètre $(t, \alpha)$. Si $Z$ est un espace de Fréchet séparable et $f : Y \rightarrow Z$ une application affine, continue et mesurable, alors $f(\eta)$ est c.l. également. Plus précisément :
\begin{itemize}
\item si $\lambda : Y \rightarrow Z$ est linéaire, $\lambda(\eta)$ est c.l. de paramètre $V \mapsto \big( t(\lambda^{-1}(V)), \alpha(\lambda^{-1}(V)) \big)$ ;
\item si $y_0 \in Y$, $\eta - y_0$ est c.l. de paramètre $V \mapsto \big( t(V) + M_V(y_0), \alpha(V) \big)$.
\end{itemize}
\end{Pro}

\subsection{Enoncé du lemme sous-additif}

Dans cette partie, on se donne un champ $\eta$ à valeurs dans $\smash{\big(Y^{\mathbb{Z}^d}, \mathcal{B}^{\otimes \mathbb{Z}^d}\big)}$. On note $\mu$ sa loi. On suppose qu'il existe un entier $\ell$ tel que $\mu$ soit invariante par les translations du sous-réseau $R_\ell = (\ell \mathbb{Z})^d$ de $\mathbb{Z}^d$. On dira que $\eta$ (ou $\mu$) est $R_\ell$-invariant. On suppose aussi que $\eta$ est c.a.d.i. et c.l. La notion de champ c.a.d.i. $R_\ell$-invariant est une généralisation naturelle de celle de suite de v.a. i.i.d. Rappelons que le théorème de Cramér, pour les moyennes de v.a. i.i.d., est conséquence du lemme sous-additif de Fekete. On montre ici que l'hypothèse de découplage asymptotique, jointe à celle de contrôle local, donne elle aussi une forme de sous-additivité, suffisante pour énoncer un résultat de type Cramér.

\medskip

Pour toute boîte $\Lambda$, on définit l'opérateur moyenne sur $\Lambda$, noté $\mathfrak{m}_\Lambda$, par
$$
\forall \xi \in Y^{\mathbb{Z}^d} \qquad \mathfrak{m}_\Lambda \xi = \frac{1}{|\Lambda|} \sum_{z \in \Lambda} \xi(z)
$$
Remarquons que, pour toute boîte $\Lambda$, $\mathfrak{m}_\Lambda$ est une application mesurable pour la tribu produit $\mathcal{B}^{\otimes \mathbb{Z}^d}$ (car la topologie sur $Y$ est à base dénombrable d'ouverts). Donnons-nous un ouvert convexe de la forme\label{convexespointes}
$$
C = y + V
$$
où $V \in \mathcal{C}_0(Y)$. Pour $\varepsilon \in ]0, 1[$, on définit
$$
C(y, \varepsilon) = y + (1 - \varepsilon)V
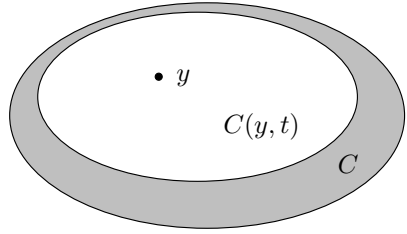
$$

\vspace{.5cm}

\begin{center}
\def\JPicScale{0.5}
\input{convpointes.pst}
\end{center}

On rappelle que, si $V$ est un convexe de $Y$, sa jauge (ou fonctionnelle de Minkowski) est l'application $M_V : Y \rightarrow [0, +\infty]$ définie par
$$
\forall y \in Y \qquad M_V(y) = \inf \{ t \geqslant 0 ; y \in tV \}
$$
De plus, si $V \in \mathcal{C}_0(Y)$, alors
$$
V = \{ y \in Y ; M_V(y) < 1 \}
$$

\begin{Th}[lemme sous-additif]\label{lsamos}
Soit $\eta$ un champ c.a.d.i., c.l. et $R_\ell$-invariant. Pour tout $\delta > 0$,
\begin{align*}
\exists M(&\varepsilon, \delta) \quad \forall m \geqslant M \quad \exists N(m, \varepsilon, \delta) \quad \forall n \geqslant N\\
&\frac{1}{|\Lambda(n)|} \log \mathbb{P} \big( \mathfrak{m}_{\Lambda(n)} \eta \in C \big) \geqslant \frac{1}{|\Lambda(m)|} \log \mathbb{P} \big( \mathfrak{m}_{\Lambda(m)} \eta \in C(y, \varepsilon) \big) - \delta
\end{align*}
\end{Th}

\subsection{Démonstration du lemme sous-additif}
La démonstration suit les étapes de celle de Pfister \cite{Pfi01} en l'adaptant à ce nouveau cadre. Si $m$ et $n \in \mathbb{N^*}$, la division euclidienne de $n$ par $m+g(m)+\ell$ s'écrit
$$
n = k \big( m+g(m)+\ell \big) + r
$$
On pave (pas complètement) $\Lambda(n)$ avec $k^d$ boîtes isométriques à $\Lambda(m+g(m)+\ell)$, que l'on note $\Lambda'_q$, pour $q \in \{ 1, \ldots , k^d \}$. Puis, pour chaque $q$, on note $\Lambda_q$ la boîte isométrique à $\Lambda(m)$ dont le plus petit élément pour l'ordre lexicographique usuel de $\mathbb{Z}^d$ --- le ``coin en bas à gauche''\,--- est le plus petit élément, pour ce même ordre, de $R_\ell \cap \Lambda_q'$. On notera que la distance entre deux boîtes $\Lambda_q$ distinctes est supérieure à $g(m)$. Soit enfin
\[
S_0 = \Lambda(n) \setminus \bigcup_{q=1}^{k^d} \Lambda_q
\]
l'ensemble des sites marginaux. Pour justifier cette dénomination, notons que
\[
\rho_{m, n} := \frac{|S_0|}{|\Lambda(n)|} \xrightarrow[\substack{m \to \infty\\n \geqslant m}]{} 0
\]
Les boîtes $\Lambda_q$ dépendent de $n$ et $m$, mais les indices seront sous-entendus.

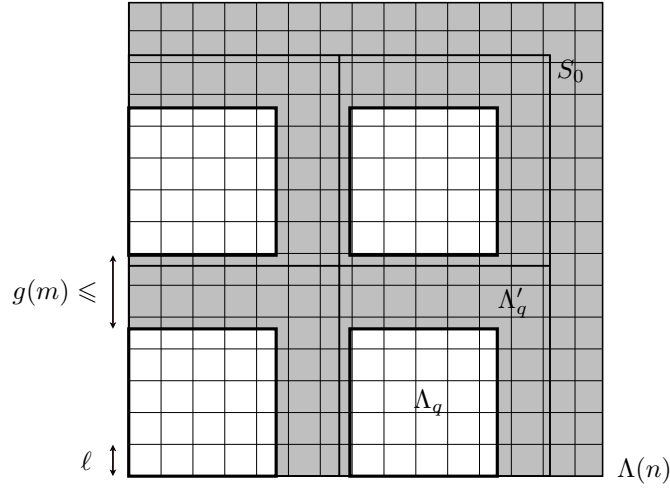
\begin{center}
\def\JPicScale{.7}
\input{lambda.pst}
\end{center}

\textbf{Remarque :} On a l'ordre de grandeur suivant :
\[
\rho_{m, n} \lesssim d \cdot \left( \frac{g(m) + \ell}{m + g(m) + \ell} + \frac{r}{n} \right)
\]
Ainsi, pour que $\rho_{m, n}$ puisse être rendu arbitrairement petit, il faut que le premier terme tende vers $0$ avec $m$, ce qui équivaut à
\[
\frac{g(m)}{m} \to 0
\]
L'hypothèse mise sur $g$ est donc utilisée ici.

\medskip

Notant $\phi$ la fonction convexe $M_V(\cdot - y)$, on a :
\begin{align*}
\phi\big(\mathfrak{m}_{\Lambda(n)} \eta\big) &\leqslant \sum_{z \in S_0} \frac{1}{|\Lambda(n)|} \phi\big( \eta(z) \big) + \sum_{q=1}^{k^d} \frac{|\Lambda(m)|}{|\Lambda(n)|} \phi\big(\mathfrak{m}_{\Lambda_q} \eta\big)\\
 &\leqslant \sum_{z \in S_0} \frac{1}{|S_0|} \rho_{m, n} \phi\big( \eta(z) \big) + \sum_{q=1}^{k^d} \frac{1}{k^d} \phi\big(\mathfrak{m}_{\Lambda_q} \eta\big)
\end{align*}
D'où
\begin{align*}
\frac{1}{|\Lambda(n)|}& \log \mathbb{P}\left(\mathfrak{m}_{\Lambda(n)} \eta \in C \right)\\
 &= \frac{1}{|\Lambda(n)|} \log \mathbb{P}\big(\phi(\mathfrak{m}_{\Lambda(n)} \eta) < 1 \big)\\
 &\geqslant \frac{1}{|\Lambda(n)|} \log \mathbb{P}\Big(\forall z \in S_0 \quad \rho_{m, n} \phi\big(\eta(z)\big) < \varepsilon \, ; \, \forall q \in \{ 1, \ldots , k^d \} \quad \phi\big(\mathfrak{m}_{\Lambda_q} \eta\big) < 1 - \varepsilon \Big)\\
 &= \frac{1}{|\Lambda(n)|} \log \mathbb{P}\Big( \forall z \in S_0 \quad \rho_{m, n} \phi\big(\eta(z)\big) < \varepsilon \, ; \, \forall q \in \{ 1, \ldots , k^d \} \quad \mathfrak{m}_{\Lambda_q} \eta \in C(y, \varepsilon) \Big)
\end{align*}
Tirons profit de l'hypothèse de contrôle local de $\eta$ : il existe $M(\varepsilon)$ et $N(m, \varepsilon)$ tels que
\[
\forall m \geqslant M(\varepsilon) \quad \forall n \geqslant N(m, \varepsilon) \quad \frac{\varepsilon}{\rho_{m, n}} > t(V)
\]
Dès lors, on obtient par récurrence sur le cardinal de $S_0$,
\begin{align*}
\mathbb{P}\Big( \forall z \in S_0& \quad \rho_{m, n} \phi\big(\eta(z)\big) < \varepsilon \, ; \, \forall q \in \{ 1, \ldots , k^d \} \quad \mathfrak{m}_{\Lambda_q} \eta \in C(y, \varepsilon) \Big)\\
 &\geqslant \alpha(V)^{|S_0|} \, \mathbb{P} \Big( \forall q \in \{ 1, \ldots , k^d \} \quad \mathfrak{m}_{\Lambda_q} \eta \in C(y, \varepsilon) \Big)
\end{align*}

Ainsi, pour $m \geqslant M(\varepsilon , \delta)$ et $n \geqslant N(m, \varepsilon , \delta)$
\begin{align}
\frac{1}{|\Lambda(n)|}& \log \mathbb{P}\left(\mathfrak{m}_{\Lambda(n)} \eta \in C \right)\nonumber\\
 &\geqslant \rho_{m, n} \log \alpha(V) + \frac{1}{|\Lambda(n)|} \log \mathbb{P}\left(\forall q \in \{ 1, \ldots , k^d \} \quad \mathfrak{m}_{\Lambda_q} \eta \in C(y, \varepsilon) \right)\label{eq1}
\end{align}

Reste à exploiter l'hypothèse de découplage asymptotique :
$$
\mathbb{P}\left(\forall q \in \{ 1, \ldots , k^d \} \quad \mathfrak{m}_{\Lambda_q} \eta \in C(y, \varepsilon) \right) \geqslant \big( e^{-c(m)} \big)^{k^d-1} \prod_{q=1}^{k^d} \mathbb{P}\big( \mathfrak{m}_{\Lambda_q} \eta \in C(y, \varepsilon)\big)
$$
Puis, utilisant le fait que $\eta$ est $R_\ell$-invariant et que $(1 - \rho_{m, n}) | \Lambda(n) | = k^d |\Lambda(m)|$, il vient
\begin{align}
\frac{1}{|\Lambda(n)|}& \log \mathbb{P}\left(\forall q \in \{ 1, \ldots , k^d \} \quad \mathfrak{m}_{\Lambda_q} \eta \in C(y, \varepsilon) \right)\nonumber\\
 &\geqslant \frac{1}{|\Lambda(m)|} \log \mathbb{P}\big( \mathfrak{m}_{\Lambda(m)} \eta \in C(y, \varepsilon)\big) - \frac{c(m)}{|\Lambda(m)|}\label{eq2}
\end{align}
Quitte à grandir $M$ et $N$, on peut faire en sorte que
\begin{equation}\label{eq3}
\forall m \geqslant M(\varepsilon, \delta) \quad \forall n \geqslant N(m, \varepsilon, \delta) \qquad \rho_{m, n} \log \alpha(V) - \frac{c(m)}{|\Lambda(m)|} \geqslant - \delta
\end{equation}
On obtient alors le résultat annoncé en combinant (\ref{eq1}), (\ref{eq2}) et (\ref{eq3}).\qed

\medskip

Intuitivement, on a d'abord fait en sorte que l'écartement entre les petites boîtes soit suffisamment petit devant la taille de ces boîtes.
Puis, on a pris suffisamment de petites boîtes pour que la proportion occupée par elles dans une grande boîte soit assez grande. Au final, on obtient une sous-additivité asymptotique entre les petites boîtes et la grande.

\medskip

\textbf{Remarque :} En combinant simplement (\ref{eq1}) et (\ref{eq2}), on obtient le résultat plus précis suivant :
\begin{align}
\exists M(\varepsilon&) \quad \forall m \geqslant M \quad \exists N(m, \varepsilon) \quad \forall n \geqslant N\nonumber\\
&\frac{1}{|\Lambda(n)|} \log \mathbb{P} \big( \mathfrak{m}_{\Lambda(n)} \eta \in C \big) \geqslant \frac{1}{|\Lambda(m)|} \log \mathbb{P} \big( \mathfrak{m}_{\Lambda(m)} \eta \in C(y, \varepsilon) \big)\nonumber\\
& \qquad\qquad\qquad\qquad\qquad\qquad\qquad - \frac{c(m)}{\Lambda(m)} + \rho_{m, n} \log \alpha(V)
\end{align}
Ici apparaissent tous les paramètres. On se souviendra que $g$ est sous-entendu dans $\rho_{m, n}$ et $t$ dans le choix de $M$ et $N$.

\section{PGD faible}\label{spgd}

On pourrait dès maintenant énoncer un principe de grandes déviations pour un champ $\eta$ c.a.d.i., c.l. et $R_\ell$-invariant, mais on se place dans le cadre plus général d'une limite projective. Cela ne complexifie pas la démonstration et s'applique à d'autres cas (notamment, $\sigma(z) = \delta_{\theta_z \omega}$, qui n'est pas c.a.d.i.).

\subsection{Limite projective}

Soient $X$ un espace vectoriel réel et
$$
\overleftarrow{\mathcal{X}} = \big( Y_i, f_i, f_{ij} \big)_{i \leqslant j}
$$
un système projectif d'espaces de Fréchet séparables, autrement dit une famille telle que

\medskip

\textsf{(PROJ$_1$)} les indices $i$ et $j$ décrivent un ensemble $(J, \leqslant)$ préordonné filtrant à droite ;

\medskip

\textsf{(PROJ$_2$)} pour tout $i \in J$, $Y_i$ est un espace de Fréchet (\emph{i.e.} un espace vectoriel localement convexe métrisable complet) séparable, $f_i$ une application linéaire de $X$ dans $Y_i$ et, pour tout $(i, j) \in J^2$ tel que $i \leqslant j$, $f_{ij}$ est une application linéaire continue de $Y_j$ dans $Y_i$ ;

\medskip

\textsf{(PROJ$_3$)} pour tout $(i, j, k) \in J^3$ tel que $i \leqslant j \leqslant k$, on a $f_{ii} = id_{Y_i}$,
$$
f_i = f_{ij} \circ f_j  \qquad \textrm{et} \qquad f_{ik} = f_{ij} \circ f_{jk}
$$

On note, pour tout $i \in J$, $\mathcal{C}_{0}(Y_i)$ un système fondamental de voisinages convexes ouverts de $0$ dans $Y_i$ et on définit
$$
\mathcal{C}_0(X) = \{ f_i^{-1}(t_i C_i) \, ; \, i \in J, \, t_i \in ]0, +\infty[, \, C_i \in \mathcal{C}_{0}(Y_i) \}
$$
On munit alors $X$ de la tribu $\mathcal{F}$ (resp. de la topologie $\tau$) engendrée par les translatés des éléments de $\mathcal{C}_0(X)$. On dit que $X$ est la \emph{limite projective} du système projectif $\overleftarrow{\mathcal{X}}$.

\medskip

%\textbf{Remarque :} La famille $\mathcal{C}_0$ satisfait aux trois axiomes \textup{\textsf{(EVLCM$_1$)}}, \textup{\textsf{(EVLCM$_2$)}} et \textup{\textsf{(EVLCM$_3$)}} du chapitre \ref{casiid}. On vérifie que, si $(X, \mathcal{C}_0, \mathcal{F}, \tau)$ est l'e.v.l.c.m. associé à $\mathcal{C}_0$, alors $\tau$ (resp. $\mathcal{F}$) est la topologie (resp. tribu) initiale pour la famille $(F_i, f_i)_{i \in J}$.

\textbf{Exemples :}
\begin{itemize}
\item si $X$ est un espace de Fréchet séparable et si $\mathcal{Y}$ est engendré par $f_0 = Id$ ($Y_0 = X$), on récupère la topologie sur $X$ ;
\item si $X$ est un espace de Fréchet séparable et si $\mathcal{Y}$ est engendré par $X'$ ($Y_i = \mathbb{R}$ pour $i \in X'$), on tombe sur la topologie $\sigma(X, X')$ ;
\item pour $X'$ dual d'un espace de Fréchet séparable, si $\mathcal{Y}$ est engendrée par $X$ ($Y_i = \mathbb{R}$ pour $i \in X$), on a la topologie $\sigma(X', X)$.
\end{itemize}

\subsection{PGD faible}

Le décor est maintenant planté. On peut énoncer\footnote{En termes de limite projective, on peut penser exploiter le théorème de Dawson-Gärtner (cf. \cite{DeZ93}) : on obtient alors un PGD pour $\sigma$, à condition que les $s_\eta$ soient de bonnes fonctions de taux, ce qui est moins bon que le résultat obtenu. Il n'est pas étonnant qu'on n'obtienne pas ainsi un résultat optimal étant donné que le théorème de Dawson-Gärtner ne se sert pas de la linéarité des $f_i$.} :

\begin{Pro}[PGD]\label{pgdmos}
Soit $\sigma$ un champ à valeurs dans $\smash{\big(X^{\mathbb{Z}^d}, \mathcal{F}^{\otimes \mathbb{Z}^d}\big)}$. On suppose que, pour tout $i \in J$, le champ $f_i(\sigma)$ est c.a.d.i., c.l. et $R_\ell$-invariant. Alors
$$
\left( \frac{1}{|\Lambda(n)|} \sum_{z \in \Lambda(n)} \sigma(z) \right)_{n \geqslant 0}
$$
vérifie un PGD faible d'entropie
$$
s(x) = \inf_{\substack{A \in \mathcal{C}\\ A \ni x}} \liminf_{n \to \infty} \frac{1}{|\Lambda(n)|} \log \mathbb{P}( m_{\Lambda(n)} \sigma \in A)
$$
\end{Pro}

%\textbf{Remarque :} En vertu de (\ref{entr}),
%\[
%s(x) = \inf_{\substack{A \in \mathcal{C}\\ A \ni x}} \underline{s}(A) = \inf_{\substack{A \in \mathcal{C}\\ A \ni x}} \overline{s}(A)
%\]

\textbf{Démonstration :} Donnons-nous $A \in \mathcal{C}$ ; il est de la forme
\[
A = x + f_i^{-1}(V)
\]
avec $V \in \mathcal{C}_0(Y_i)$. Notons $y = f_i(x)$ et définissons
\[
C = y + V
\]
puis, pour $\varepsilon \in ]0, 1[$,
\[
A(x, \varepsilon) = x + (1 - \varepsilon) U \quad \textrm{et} \quad C(y, \varepsilon) = y + (1 - \varepsilon) V
\]

La famille $\eta = f_i(\sigma)$ vérifie les hypothèses du lemme sous-additif \ref{lsamos} : on en déduit que, pour tout $\delta > 0$,
\begin{align*}
\exists M(&\varepsilon, \delta) \quad \forall m \geqslant M \quad \exists N(m, \varepsilon, \delta) \quad \forall n \geqslant N\\
&\frac{1}{|\Lambda(n)|} \log \mathbb{P} \big( \mathfrak{m}_{\Lambda(n)} \eta \in C \big) \geqslant \frac{1}{|\Lambda(m)|} \log \mathbb{P} \big( \mathfrak{m}_{\Lambda(m)} \eta \in C(y, \varepsilon) \big) - \delta
\end{align*}
Puis, en faisant tendre $n$, puis $m$, vers l'infini, et enfin $\delta$ vers $0^+$,
\[
\liminf_{n \to \infty} \frac{1}{|\Lambda(n)|} \log \mathbb{P} \big( \mathfrak{m}_{\Lambda(n)} \eta \in C \big) \geqslant \limsup_{m \to \infty} \frac{1}{|\Lambda(m)|} \log \mathbb{P} \big( \mathfrak{m}_{\Lambda(m)} \eta \in C(y, \varepsilon) \big)
\]

Pour toute boîte finie $\Lambda$, la linéarité de $f_i$ assure
\[
\{ \mathfrak{m}_\Lambda \sigma \in A \} = \{ f_i(\mathfrak{m}_\Lambda \sigma) \in C \} = \{ \mathfrak{m}_\Lambda \eta \in C \}
\]
et les égalités analogues avec $A(x, \varepsilon)$ et $C(y, \varepsilon)$. On en déduit
\[
\underline{s}(A) \geqslant \overline{s}(A(x, \varepsilon))
\]
La conclusion est un résultat classique de grandes déviations, qu'on trouvera par exemple en \cite[Proposition 3.5.]{Pfi01}.\qed

\begin{Pro}
L'entropie $s$ est s.c.s. et concave.
\end{Pro}
\textbf{Démonstration :} Elle se fait comme dans le cas indépendant. La semi-continuité est immédiate : si $u \in \mathbb{R}$ et $x \in \{ s < u \}$, il existe $A \in \mathcal{C}$ contenant $x$ tel que $\underline{s}(A) < u$. Alors, pour tout $y \in A$, $s(y) \leqslant \underline{s}(A) < u$, donc $x \in A \subset \{ s < u \}$ et $\{ s < u \}$ est ouvert. Pour la concavité, il suffit alors de prouver que
\[
s\left( \frac{x+y}{2} \right) \geqslant \frac{s(x)+s(y)}{2}
\]
La semi-continuité permet alors de passer des inégalités
\[
s((1 - u)x + uy) \geqslant (1 - u)s(x)+ us(y)
\]
pour $u$ dyadique à celles pour $u \in [0, 1]$. La preuve est un raffinement de la démonstration du lemme sous-additif (\ref{lsamos}). Soit $A \in \mathcal{C}$ contenant $\frac{1}{2}(x+y)$. Par continuité des applications d'espace vectoriel, pour tout $\varepsilon \in ]0, 1[$, il existe $A_x$ et $A_y \in \mathcal{C}$, contenant respectivement $x$ et $y$, tels que
\[
\frac{1}{2}A_x + \frac{1}{2}A_y \subset A\left(\frac{1}{2}(x+y), \varepsilon\right)
\]
Puis
%On se restreint, dans le calcul suivant, aux couples $(m, n)$ tels que $k$ est pair.
\begin{align*}
\frac{1}{|\Lambda(n)|} \log &\, \mathbb{P}\left(\mathfrak{m}_{\Lambda(n)} \sigma \in A \right)\\
 &\geqslant \frac{1}{|\Lambda(n)|} \log \mathbb{P}\Bigg(\forall z \in \Lambda_0 \quad \rho_{m, n} \phi(\eta(z)) < \varepsilon \, ;\\
&\qquad\qquad\qquad\qquad\bigcap_{q \textrm{ pair}} \{ \mathfrak{m}_{\Lambda_q} \sigma \in A_x \} \cap \bigcap_{q \textrm{ impair}} \{ \mathfrak{m}_{\Lambda_q} \sigma \in A_y \} \Bigg)
\end{align*}
En poursuivant comme précédemment, on obtient
\[
\underline{s}(A) \geqslant \frac{1}{2} \big( \underline{s}(A_x) + \underline{s}(A_y) \big) \geqslant \frac{1}{2} \big( s(x) + s(y) \big)
\]
puis l'inégalité désirée en passant à l'infimum sur les $A \ni \frac{1}{2}(x+y)$.

\section{Pression}\label{spression}

Soient $\overleftarrow{\mathcal{X}} = (Y_i, f_i, f_{ij})_{i \leqslant j}$ un système projectif d'espaces de Fréchet séparables et $X$ sa limite projective. Soit $\sigma$ un champ à valeurs dans $\smash{\big(X^{\mathbb{Z}^d}, \mathcal{F}^{\otimes \mathbb{Z}^d}\big)}$ tel que, pour tout $i \in J$, le champ $f_i(\sigma)$ soit c.a.d.i., c.l. et $R_\ell$-invariant. On définit, pour toute boîte $\Lambda \subset \mathbb{Z}^d$,
\[
\forall \lambda \in X^* \quad p_\Lambda(\lambda) = \frac{1}{|\Lambda|} \log \mathbb{E} \left[ \exp \left\langle \lambda \bigg| \sum_{z \in \Lambda} \sigma(z) \right\rangle \right]
\]

Notons $X^*$ le dual topologique de $X$ Soit $\lambda \in X^*$. Il existe un voisinage convexe de $0$ sur lequel $|\lambda| < 1$, voisinage de la forme $f_i^{-1}(C_i)$. Par linéarité, on remarque que
\[
\forall x, x' \in X \quad f_i(x) = f_i(x') \Rightarrow \lambda(x) = \lambda(x')
\]
Ainsi, on peut définir la section
\[
\lambda_i \in Y_i^* : y \mapsto \lambda\big(f_i^{-1}(y)\big)
\]
de sorte que $\lambda = \lambda_i \circ f_i$. Donc
\[
X^* = \{ \lambda_i \circ f_i ; i \in I, \lambda_i \in Y_i^* \}
\]

En particulier, pour tout $\lambda \in X^*$, écrivant $\lambda$ sous la forme $\lambda_i \circ f_i$, on voit que le champ $\eta = \langle \lambda | \sigma \rangle = \langle \lambda_i | f_i(\sigma) \rangle$ est c.a.d.i., c.l. (cf. \ref{scadimos} et \ref{sclmos}) et $R_\ell$-invariant.

\subsection{Pression en volume infini}

On veut maintenant définir la pression $p$ en volume infini. Tout se passe à merveille ainsi que le montre la proposition suivante :
\begin{Pro}
Pour tout $\lambda \in X^*$, la suite $\big(p_{\Lambda(n)}(\lambda)\big)_{n \in \mathbb{N}}$ converge dans $[- \infty , + \infty]$. On note $p(\lambda)$ sa limite.
\end{Pro}

\textbf{Démonstration :} La preuve est analogue à celle du lemme sous-additif \ref{lsamos}. On définit $\eta(z) = \langle \lambda | \sigma(z) \rangle$ ; on note $(g, c)$ (resp. $(t, \alpha)$) le paramètre de découplage (resp. de contrôle local) de $\eta$. On écrit alors
\[
p_{\Lambda(n)}(\lambda) = \frac{1}{|\Lambda(n)|} \log \mathbb{E}\big[ P_{m, n} \cdot R_{m, n}\big]
\]
avec un terme principal
\[
P_{m, n} = \prod_{q=1}^{k^d} \exp \bigg( \sum_{z \in \Lambda_q} \eta(z) \bigg)
\]
et un terme résiduel
\[
R_{m, n} = \mathbb{E}\left[\exp \bigg(\sum_{z \in S_0} \eta(z)\bigg) \Bigg| \prod_{q=1}^{k^d} \exp \bigg( \sum_{z \in \Lambda_q} \eta(z) \bigg)\right]
\]

Pour qu'il y ait convergence, il suffit de s'assurer que
\begin{itemize}
\item l'on a une minoration du terme résiduel de la forme
\[
R_{m, n} \geqslant e^{- \delta |\Lambda(n)|}
\]
pour $m$ et $n$ assez grands.
\item l'hypothèse de découplage s'applique bien pour traiter le terme principal.
\end{itemize}

\textbf{Terme résiduel :} Encore une fois, on va utiliser le contrôle local.

\begin{Le}
Il existe $\beta > 0$ tel que, pour toute boîte finie $\Lambda \subset \mathbb{Z}^d \setminus \{ 0 \}$,
\[
\mathbb{E}\big[ e^{\eta(0)} \big| \eta_\Lambda \big] \geqslant \beta
\]
\end{Le}
\textbf{Démonstration du lemme :} On note que, pour tous $\Lambda$ finie et $t > 0$,
\begin{align*}
\mathbb{E}\big[ e^{\eta(0)} \big| \eta_\Lambda \big] &\geqslant \mathbb{E}\big[ e^{\eta(0)} 1_{|\eta(0)| < t} \big| \eta_\Lambda \big]\\
 &\geqslant e^{-t} \mathbb{P}\big( |\eta(0)| < t \big| \eta_\Lambda \big)
\end{align*}
En utilisant l'hypothèse de contrôle local pour $\eta$ et $V_1 = ]-1, 1[$, on note $t = t(V_1)$ et $\alpha = \alpha(V_1)$, de sorte que
\begin{align*}
\mathbb{P}_{\eta_\Lambda}\!-\!p.s. \qquad \mathbb{P}\big( |\eta(0)| < t \big| \eta_\Lambda = y \big) &= \lim_{r \to 0^+} \mathbb{P}\Big( |\eta(0)| < t \Big| \eta_\Lambda \in B(y, r) \Big)\\
 &\geqslant \alpha
\end{align*}
L'égalité de la première ligne est due au fait que l'ensemble de Lebesgue d'une fonction de $L^1(\mathbb{R}^{\Lambda} ; \mathbb{R})$ est de mesure pleine (cf. par exemple \cite[théorème 8.8]{Rud66}). Ainsi, $\beta = e^{-t} \alpha > 0$ convient.\qed

\medskip

On déduit du lemme, par récurrence, que
\[
R_{m, n} \geqslant \beta^{|S_0|} \geqslant e^{-\delta |\Lambda(n)|}
\]
pour $m \geqslant M(\delta)$ et $n \geqslant N(m, \delta)$.

\medskip

\textbf{Terme principal :} On a
\[
\left\{ \exp \sum_{z \in \Lambda_q} \eta(z) > t \right\} \in \mathcal{C}\big(\mathbb{R}^{|\Lambda_q|}\big)
\]
(il s'agit en fait d'un demi-espace ouvert). Ainsi,
\begin{align*}
\frac{1}{|\Lambda(n)|} \log &\, \mathbb{E}\left[ \prod_{q=1}^{k^d} \exp \bigg( \sum_{z \in \Lambda_q} \eta(z) \bigg)\right]\\
 &= \frac{1}{|\Lambda(n)|} \log \mathbb{E}\left[ \prod_{q=1}^{k^d} \int_{t_q \geqslant 0} 1_{\exp \Big( \sum\limits_{z \in \Lambda_q} \eta(z) \Big) > t_q}\right]\\
 &= \frac{1}{|\Lambda(n)|} \log \int_{t_1, \ldots , t_{k^d} \geqslant 0}\mathbb{E}\left[ \prod_{q=1}^{k^d} 1_{\exp \Big( \sum\limits_{z \in \Lambda_q} \eta(z) \Big) > t_q}\right]\\
 &\geqslant \frac{1}{|\Lambda(n)|} \log \int_{t_1, \ldots , t_{k^d} \geqslant 0} \big(e^{-c(m)}\big)^{k^d-1} \prod_{q=1}^{k^d} \mathbb{E}\left[1_{\exp \Big( \sum\limits_{z \in \Lambda_q} \eta(z) \Big) > t_q}\right]\\
 &\geqslant \frac{k^d}{|\Lambda(n)|} \log \mathbb{E}\left[\exp \bigg( \sum\limits_{z \in \Lambda(m)} \eta(z) \bigg)\right] - \frac{c(m)}{|\Lambda(m)|}
\end{align*}

On obtient finalement que, quitte à grandir $M$ et $N$ (pour absorber le terme $c(m)/|\Lambda(m)|$),
\begin{align}
\forall m \geqslant M(\delta) \quad \forall n \geqslant N(m&, \delta)\nonumber\\
p_{\Lambda(n)}(\lambda)& \geqslant \frac{1}{|\Lambda(n)|} \log \mathbb{E}\Big[ P_{m, n} \cdot e^{-\delta |\Lambda(n)|}\Big]\nonumber\\
 &\geqslant (1 - \rho_{m, n}) p_{\Lambda(m)}(\lambda) - \delta
\end{align}

D'où, en prenant la limite inférieure en $n$, puis la limite supérieure en $m$ et, enfin, la limite quand $\delta \to 0^+$,
\[
\liminf_{n \to \infty} p_{\Lambda(n)}(\lambda) \geqslant \limsup_{m \to \infty} p_{\Lambda(m)}(\lambda)
\]
autrement dit $p_{\Lambda(n)}(\lambda)$ converge.\qed

\medskip

\textbf{Remarque :} Voici, comme pour $s$, la version de l'inégalité sous-additive avec tous les paramètres :
\begin{align}\label{plsa}
\forall n \geqslant N(m)\nonumber\\
p_{\Lambda(n)}(\lambda) &\geqslant (1 - \rho_{m, n}) p_{\Lambda(m)}(\lambda) - \frac{c(m)}{|\Lambda(m)|} - \rho_{m, n}\big(t(V_1) - \log \alpha(V_1)\big)
\end{align}
On se souviendra que les paramètres intervenant ici sont ceux de $\langle \lambda | \sigma \rangle$, et dépendent donc de $\lambda$.

\begin{Pro}[Convexité de $p$]\label{pconv}
La pression $p$ est convexe.
\end{Pro}
\textbf{Démonstration :} L'inégalité de Hölder assure que les $p_\Lambda$ sont convexes. Puis, sachant qu'une limite supérieure de fonctions convexes est convexe, $p$ l'est également.

\section{Pression et entropie}

\subsection{Inégalité de Tchebychev}

L'inégalité de Tchebychev exponentielle (cf. \cite[lemme 12.6]{Cer07}) prend ici la forme suivante :
\begin{Pro}[Inégalité de Tchebytcheff exponentielle]\label{ite}
Pour toute boîte $\Lambda \subset \mathbb{Z}^d$ et $A \in \mathcal{B}$,
\[
\mathbb{P}(\mathfrak{m}_\Lambda \sigma \in A) \leqslant \exp \left[ - |\Lambda| \sup_{\lambda \in X^*} \Big( \inf_{x \in A} \langle \lambda | x \rangle - p_\Lambda(\lambda) \Big) \right]
\]
\end{Pro}

On en déduit, en considérant $\Lambda = \Lambda(n)$ et en faisant tendre $n$ vers $\infty$ :
\[
\underline{s}(A) \leqslant \liminf \inf_{\lambda \in X^*} \Big( p_{\Lambda(n)}(\lambda) - \inf_{a \in A} \langle \lambda | a \rangle \Big) \leqslant \inf_{\lambda \in X^*} \Big( p(\lambda) - \inf_{a \in A} \langle \lambda | a \rangle \Big)
\]

Puis,
\begin{align*}
s(x) &= \inf_{\substack{A \in \mathcal{C}\\ A \ni x}} \underline{s}(A)\\
 &\leqslant \inf_{\substack{A \in \mathcal{C}_w\\ A \ni x}} \underline{s}(A) = \underline{s}_m(x)\\
 &\leqslant \inf_{\substack{A \in \mathcal{C}_w\\ A \ni x}} \inf_{\lambda \in X^*} \Big( p(\lambda) - \inf_{a \in A} \langle \lambda | a \rangle \Big)\\
 &= \inf_{\lambda \in X^*} \inf_{\substack{A \in \mathcal{C}_w\\ A \ni x}} \Big( p(\lambda) - \inf_{a \in A} \langle \lambda | a \rangle \Big)\\
 &= - \sup_{\lambda \in X^*} \Big( \langle \lambda | x \rangle - p(\lambda) \Big)\\
 &= - p^*(x)
\end{align*}

Résumons :
\begin{Th}
On a toujours l'inégalité $s \leqslant - p^*$.
\end{Th}

On voudrait maintenant savoir dans quels cas l'égalité est vérifiée.

\subsection{\'Egalité $s = -p^*$ : cas borné}

\begin{Th}\label{cspetoile}
Si $\sigma(0)$ est à valeurs dans un ensemble faiblement borné, on a l'égalité
\[
s(\sigma ; x) = - p^*(\sigma ; x)
\]
\end{Th}

\textbf{Démonstration} : L'égalité $p = (-s)^*$ est conséquence d'une version du lemme de Varadhan en remarquant que
\[
p_{\Lambda(n)} = \frac{1}{|\Lambda(n)|} \log \int e^{|\Lambda(n)| \langle \lambda | x \rangle} d \nu_n(x)
\]
où $\nu_n$ est la loi de $m_{\Lambda(n)} \sigma$. L'égalité duale découle alors du fait que que $s$ est concave et semi-continue supérieurement (cf. par exemple \cite[Proposition 5.e]{Mor67}).\qed

\subsection{\'Egalité $s = -p^*$ : cas c.a.d.i. et c.l.}

\begin{Th}\label{lspetoile}
Soit $(Y, \tau)$ un espace de Fréchet séparable et $\mathcal{B}$ sa tribu borélienne. Soit $\eta$ un champ à valeurs dans $\smash{\big(Y^{\mathbb{Z}^d}, \mathcal{B}^{\otimes \mathbb{Z}^d}\big)}$ c.a.d.i., c.l. et $R_\ell$-invariant. Alors,
\[
s(\eta ; y) = - p^*(\eta ; y)
\]
\end{Th}

\textbf{Démonstration :} L'idée de la démonstration est de se ramener au théorème \ref{cspetoile} en conditionnant les sites de $\eta$ à être dans des compacts. Comme il y a une infinité de sites, on est obligé de prendre des précautions pour ne pas manipuler des conditionnements dégénérés. C'est l'objet de la définition suivante :
\begin{De}
Si $\mu$ est une loi sur $Y^{\mathbb{Z}^d}$ et $m \in \mathbb{N}$, on note $\mu_{\Lambda(m)}$ la loi sur $Y^{\mathbb{Z}^d}$ définie par
\[
\mu_{\Lambda(m)} := \big( \mu|_{\Lambda(m)} \big)^{\otimes \mathbb{Z}^d}
\]
où $\mu|_{\Lambda(m)}$ désigne la marginale de $\mu$ sur $Y^{\Lambda(m)}$. On définit aussi, pour $K$ compact de $Y$, $\mu_{\Lambda(m)}^K$ la loi sur $Y^{\mathbb{Z}^d}$ donnée par
\[
\mu_{\Lambda(m)}^K := \Big( \mu|_{\Lambda(m)}\big(\cdot | \forall z \in \Lambda(m) \quad \xi(z) \in K \big) \Big)^{\otimes \mathbb{Z}^d}
\]
(quand le conditionnement a un sens).
\end{De}

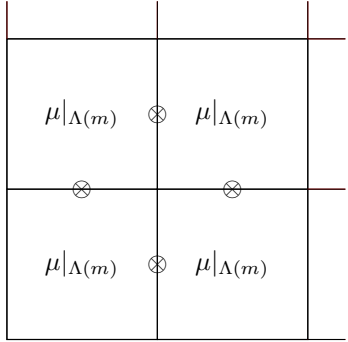
\begin{center}
\def\JPicScale{.5}
\input{mulambda.pst}
\end{center}

\textbf{Remarque :} Comme nous utiliserons ces notations, nous préférerons l'usage de $\mu$ (resp. $\xi$) à celui de $\mathbb{P}$ (resp. $\eta$).

\medskip

Soit $\mu$ la loi de $\eta$, $(g, c)$ et $(t, \alpha)$ ses paramètres respectifs de découplage et de contrôle local. On reprend les notations du lemme sous-additif (cf. p. \pageref{convexespointes}) :
\[
C = y + V \quad \textrm{et} \quad C(y, \varepsilon) = y + (1 - \varepsilon)V
\]
où $V \in \mathcal{C}_0(Y)$ et $\varepsilon \in ]0, 1[$. On note aussi $\varphi = M_V(\cdot - y)$. La démonstration sera faite en deux étapes :

\medskip

\textbf{\'Etape 1 :} Si $(K_m)_{m \in \mathbb{N}}$ est une suite de compacts vérifiant
\begin{equation}\label{compexh}
\lim_{m \to \infty} \mu\big( K_m^{\Lambda(m+g(m)+\ell)} \big) = 1
\end{equation}
alors
\begin{align*}
\forall y \in &\, Y \quad \forall \varepsilon \in ]0, 1[\\
&s(\mu ; y) \geqslant \inf_{\substack{C \in \mathcal{C}(Y)\\ C \ni y}} \limsup_{m\to\infty} \sup_{x \in C(y, \varepsilon)} \big[ - p^*\big( \mu_{\Lambda(m+g(m)+\ell)}^{K_m} ; x \big) \big]
\end{align*}

\textbf{\'Etape 2 :} Grâce au théorème de Mosco, on exhibe une suite de compacts $(K_m)_{m \in \mathbb{N}}$ vérifiant (\ref{compexh}) et telle que 
\[
\liminf_{m\to\infty} \sup_{x \in C(y, \varepsilon)} \big[ - p^*\big( \mu_{\Lambda(m+g(m)+\ell)}^{K_m} ; x \big) \big] \geqslant \sup_{x \in C(y, \varepsilon)} \big[ - p^*(\mu ; x) \big]
\]

On conclut alors ainsi :
\[
s(\mu ; y) \geqslant \inf_{\substack{C \in \mathcal{C}(Y)\\ C \ni y}} \sup_{x \in C(y, \varepsilon)} \big[ - p^*(\mu ; x) \big] \geqslant - p^*(\mu ; y)
\]

\textbf{\'Etape 1 :} Par un raisonnement analogue à la démonstration du théorème \ref{lsamos}, on montre que
\begin{align*}
\frac{1}{|\Lambda(n)|}& \log \mu \big( \mathfrak{m}_{\Lambda(n)} \xi \in C \big)\\ 
 &\geqslant \frac{1}{|\Lambda(n)|} \log \mu\Big( (1 - \rho_{m, n}) \varphi \big(\mathfrak{m}_{\bigcup\Lambda_q} \xi \big) < 1 - \textstyle{\frac{\varepsilon}{2}} \Big) + \rho_{m, n} \log \alpha(V)\\
 &\geqslant \frac{1}{|\Lambda(n)|} \log \mu_{\Lambda(m+g(m)+\ell)} \Big( (1 - \rho_{m, n}) \varphi \big(\mathfrak{m}_{\bigcup\Lambda_q} \xi \big) < 1 - \textstyle{\frac{\varepsilon}{2}} \Big)\\
 & \qquad\qquad\qquad\qquad\qquad\qquad\qquad\qquad - \frac{c(m)}{\Lambda(m)} + \rho_{m, n} \log \alpha(V)
\end{align*}
En effet, utilisant le fait que $Y$ est à base dénombrable d'ouverts et adaptant les idées de la démonstration de \cite[lemme 11.2.]{Cer07}, on peut écrire l'ensemble
$$
\left\{ (y_1, \ldots , y_{k^d}) \in Y^{k^d} \, ; \, \frac{y_1 + \cdots + y_{k^d}}{k^d} \in \big( 1 - \textstyle{\frac{\varepsilon}{2}}\big) V  \right\}
$$
comme union dénombrable de produits d'ouverts convexes
$$
\bigsqcup_{i \in \mathbb{N}} V_1^i \times \cdots \times V_{k^d}^i
$$
et on en déduit :
\begin{align*}
\mu\Big( (1 - &\rho_{m, n}) \varphi \big(\mathfrak{m}_{\bigcup\Lambda_q} \xi \big) < 1 - \textstyle{\frac{\varepsilon}{2}} \Big)\\
 &= \mu\Big( (1 - \rho_{m, n}) \big(\mathfrak{m}_{\bigcup\Lambda_q} \xi - y \big) \in (1 - \textstyle{\frac{\varepsilon}{2}}) V \Big)\\
 &= \mu\left( \bigsqcup_{i \in \mathbb{N}} \prod_{q=1}^{k^d} \big\{ (1 - \rho_{m, n}) \big(\mathfrak{m}_{\Lambda_q} \xi - y \big) \in V_q^i \big\} \right)\\
 &\geqslant e^{- (k^d - 1) c(m)} \mu_{\Lambda(m+g(m)+\ell)}\left( \bigsqcup_{i \in \mathbb{N}} \prod_{q=1}^{k^d} \big\{ (1 - \rho_{m, n}) \big(\mathfrak{m}_{\Lambda_q} \xi - y \big) \in V_q^i \big\} \right)\\
 &= e^{- (k^d - 1) c(m)} \mu_{\Lambda(m+g(m)+\ell)}\Big( (1 - \rho_{m, n}) \big(\mathfrak{m}_{\bigcup\Lambda_q} \xi - y \big) \in (1 - \textstyle{\frac{\varepsilon}{2}}) V \Big)\\
 &= e^{- (k^d - 1) c(m)} \mu_{\Lambda(m+g(m)+\ell)}\Big( (1 - \rho_{m, n}) \varphi \big(\mathfrak{m}_{\bigcup\Lambda_q} \xi \big) < 1 - \textstyle{\frac{\varepsilon}{2}} \Big)
\end{align*}

Enfin, soit $(K_m)_{m\in\mathbb{N}}$ une suite de compacts mesurables de $Y$. On a
\begin{align*}
\frac{1}{|\Lambda(n)|}\log & \,\,\mu_{\Lambda(m+g(m)+\ell)}\Big( (1 - \rho_{m, n}) \varphi \big(\mathfrak{m}_{\bigcup\Lambda_q} \xi \big) < 1 - \textstyle{\frac{\varepsilon}{2}} \Big)\\
 &\geqslant \frac{1}{|\Lambda(n)|} \log \mu_{\Lambda(m+g(m)+\ell)}^{K_m} \Big( \mathfrak{m}_{\Lambda(n)} \xi \in C(y, \varepsilon) \, ; \, \rho_{m, n} \varphi \big( - \mathfrak{m}_{\Lambda_0} \xi \big) < \textstyle{\frac{\varepsilon}{2}} \Big)\\
 &\quad + \frac{1}{|\Lambda(n)|} \log \mu_{\Lambda(m+g(m)+\ell)} \big( \forall z \in \Lambda(n) \quad \xi(z) \in K_m \big)
\end{align*}

Comme $-y - K_m$ est absorbé par $V$,
\begin{align*}
\liminf_{n \to \infty} \frac{1}{|\Lambda(n)|} \log&\, \mu_{\Lambda(m+g(m)+\ell)}^{K_m} \Big( \mathfrak{m}_{\Lambda(n)} \xi \in C(y, \varepsilon) \, ; \, \rho_{m, n} \varphi \big( - \mathfrak{m}_{\Lambda_0} \xi \big) < \textstyle{\frac{\varepsilon}{2}} \Big)\\
 &= \liminf_{n \to \infty}\frac{1}{|\Lambda(n)|} \log \mu_{\Lambda(m+g(m)+\ell)}^{K_m} \Big( \mathfrak{m}_{\Lambda(n)} \xi \in C(y, \varepsilon) \Big)
\end{align*}

D'autre part,
\begin{align*}
\liminf_{n \to \infty}& \frac{1}{|\Lambda(n)|} \log \mu_{\Lambda(m+g(m)+\ell)} \big( \forall z \in \Lambda(n) \quad \xi(z) \in K_m \big)\\
 &\geqslant \liminf_{n \to \infty} \frac{1}{|\Lambda(n)|} \log \mu_{\Lambda(m+g(m)+\ell)} \big( \forall z \in \Lambda\big((k+1)\cdot(m+g(m)+\ell)\big) \quad \xi(z) \in K_m \big)\\
 &\geqslant \frac{1}{|\Lambda(m+g(m)+\ell)|} \log \mu \big( \forall z \in \Lambda(m+g(m)+\ell) \quad \xi(z) \in K_m \big)
\end{align*}

Pour tout $z \in \mathbb{Z}^d$, la loi $\nu$ de $\eta(z)$ est tendue, donc il existe une suite croissante $(K_m)_{m \geqslant 1}$ de compacts de $Y$ telle que, pour tout $m \geqslant 1$,
$$
\nu(K_m) \geqslant 1 - \frac{1}{m |\Lambda(m+g(m)+\ell)|}
$$
Alors,
\[
\lim_{m\to\infty} \mu\big(K_m^{\Lambda(m+g(m)+\ell)}\big) = 1
\]
Introduisons la notation
\[
\mu_m := \mu_{\Lambda(m+g(m)+\ell)}^{K_m}
\]
On obtient finalement, en faisant tendre $n$, puis $m$, vers l'infini :
\begin{align*}
\forall \varepsilon > 0 \qquad \underline{s}(\mu ; C)
 & \geqslant \limsup_{m \to \infty} \, \underline{s}\big(\mu_m ; C(y, \varepsilon) \big)\\
 & \geqslant \limsup_{m \to \infty} \, \sup_{x \in C(y, \varepsilon)}s\big(\mu_m ; x \big)\\
 &= \limsup_{m \to \infty} \sup_{x \in C(y, \varepsilon)} \big[ - p^*\big( \mu_m ; x \big)\big]
\end{align*}
L'égalité finale est conséquence du théorème \ref{cspetoile}. On conclut cette première étape en passant à l'infimum sur $C \ni y$.\qed

\medskip

%%%%%%%%%%%%
%[La condition n'est pas nécessaire : que donne l'hypothèse minimale
%\[
%\frac{1}{|\Lambda(m+g(m))|} \log \mu \big( \forall z \in \Lambda(m+g(m)) \quad x_z \in K_m \big) \to 0 \quad ?
%\]
%C'est une condition moins contraignante que la tension exponentielle. En effet, on l'obtient par procédé diagonal dans l'hypothèse de tension exponentielle. On peut s'intéresser à la tension exponentielle qui donne peut-être le PGD fort et la coercivité de $s$ (comme dans le cas indépendant) : que donne alors $s = -p^*$ via l'approche $cosupp(\mu)$ faiblement borné au lieu de compact ?]
%%%%%%%%%%%%%

\textbf{\'Etape 2 :} Montrons maintenant que
\begin{equation}\label{pmosco}
\limsup_{m \to \infty} \inf_{x \in C(y, \varepsilon)}  p^*\big( \mu_m ; x \big) \leqslant \inf_{x \in C(y, \varepsilon)} p^*( \mu ; x)
\end{equation}

Nous allons voir que cette inégalité est conséquence du théorème suivant, dont on trouvera une démonstration dans \cite{Zab92b}.
\begin{Th}[Mosco]\label{mosco}
Soit $Y$ un e.v.l.c. séparable et métrisable. On se donne $f~:~[-\infty,+\infty]~\rightarrow~Y^*$ et $(f_m)_{m \in \mathbb{N}} \in [-\infty, +\infty]^{Y^*}$ une suite de fonctions convexes s.c.i., uniformément propre. On suppose que $f_m \xrightarrow[m\to\infty]{M_2} f$, \emph{i.e.}
\[
\forall \lambda \in Y^* \quad \forall \lambda_m \xrightarrow[m\to\infty]{\sigma(Y^*, Y)} \lambda \qquad \liminf_{m\to\infty} f_m(\lambda_m) \geqslant f(\lambda)
\]
Alors $f_m^* \xrightarrow[m\to\infty]{M_1} f^*$, \emph{i.e.}
\[
\forall y \in Y \quad \exists y_m \xrightarrow[m\to\infty]{\tau(Y^*, Y)} y \qquad \limsup_{m\to\infty} f_m^*(y_m) \leqslant f^*(y)
\]
\end{Th}

\textbf{Conclusion de l'étape 2} (en admettant que $f = p(\mu ; \cdot)$ et $f_m = p(\mu_m ; \cdot)$ vérifient les hypothèses du théorème \ref{mosco}) \textbf{:} Soit $x_0 \in C(y, \varepsilon)$. Il existe une suite $y_m \to x_0$ telle que
\[
\limsup_{m\to\infty} p^*(\mu_m ; y_m) \leqslant p^*(\mu ; x_0)
\]
Or, pour tout $m$ assez grand, $y_m \in C(y, \varepsilon)$, donc
\[
\exists m_0 \quad \forall m \geqslant m_0 \qquad \inf_{x \in C(y, \varepsilon)} p^*(\mu_m ; x) \leqslant p^*(\mu_m ; y_m)
\]
D'où
\[
\limsup_{m\to\infty} \inf_{x \in C(y, \varepsilon)} p^*(\mu_m ; x) \leqslant \limsup_{m\to\infty} p^*(\mu_m ; y_m) \leqslant p^*(\mu ; x_0)
\]
On obtient (\ref{pmosco}) en passant à l'infimum sur $x_0 \in C(y, \varepsilon)$.\qed

\medskip

Vérifions donc que $f = p(\mu ; \cdot)$ et $f_m = p( \mu_m ; \cdot)$ satisfont les hypothèses du théorème de Mosco \ref{mosco}. Les fonctions $f_m$ sont convexes d'après la proposition \ref{pconv}. Pour voir qu'elles sont s.c.i. (pour la topologie faible, par exemple, ou la topologie sur $X$ ; cela ne change rien puisqu'elles sont convexes), il suffit de constater (cf. la proposition \ref{pmulambda} ci-dessous) que
\[
p\big( \mu_{\Lambda(m+g(m)+\ell)}^{K_m} ; \lambda \big) = p_{\Lambda(m+g(m)+\ell)}\big(\mu_{\Lambda(m+g(m)+\ell)}^{K_m} ; \lambda\big)
\]
ce qui permet de se ramener à \cite[lemme 12.1.]{Cer07}.

\begin{Pro}\label{pmulambda}
Pour tout $j \in \mathbb{N}$, on a :
\[
p(\mu_{\Lambda(j)} ; \lambda) = p_{\Lambda(j)}(\mu_{\Lambda(j)} ; \lambda)
\]
\end{Pro}
\textbf{Démonstration :} En effet, on peut écrire :
\begin{align*}
p(\mu_{\Lambda(j)} ; \lambda) &= \lim_{n \to \infty} p_{\Lambda(n)}(\mu_{\Lambda(j)} ; \lambda)\\
 &= \lim_{n=kj \to \infty} p_{\Lambda(n)}(\mu_{\Lambda(j)} ; \lambda)\\
 &= \lim_{n=kj \to \infty} \frac{1}{|\Lambda(n)|} \log \int \exp \left\langle \lambda \bigg| \sum_{z \in \Lambda(n)} \xi(z) \right\rangle d\mu_{\Lambda(j)}(\xi)\\
 &= \lim_{n=kj \to \infty} \frac{1}{|\Lambda(n)|} \log \prod_{q=1}^{k^d} \int \exp\left\langle \lambda \bigg| \sum_{z \in \Lambda'_q} \xi(z) \right\rangle d \mu_{\Lambda(j)}(\xi)\\
 &= p_{\Lambda(j)}(\mu_{\Lambda(j)} ; \lambda)
\end{align*}
où les $\Lambda'_q$ sont les boîtes translatées de $\Lambda(j)$ qui pavent naturellement $\Lambda(kj) = \Lambda(n)$.\qed

\begin{De}
Soit $(f_m)_{m \in \mathbb{N}}$ une famille de fonctions convexes propres (\emph{i.e.} pour tout $m$, $f_m > - \infty$ et il existe $\lambda_m \in Y^*$ tel que $f_m(\lambda_m) < \infty$). On dit que $(f_m)_{m \in \mathbb{N}}$ est uniformément propre s'il existe une suite $(\lambda_m)_{m \in \mathbb{N}}$ $\sigma(Y^*, Y)$-relativement compacte telle que
\[
\sup_{m\in\mathbb{N}} f_m(\lambda_m) < \infty
\]
\end{De}

Il est alors immédiat que la famille de fonctions $\big(p(\mu_m ; \cdot)\big)_{m \in \mathbb{N}}$ est uniformément propre, car leur valeur en $0$ est $0$ ($\lambda_m = 0$ convient). Reste à montrer la convergence $(M_2)$. Pour ce faire, soient $\lambda \in Y^*$ et $(\lambda_m)_{m\in\mathbb{N}}$ tels que $\lambda_m$ converge vers $\lambda$ pour la topologie $\sigma(Y^*, Y)$. On procède en deux sous-étapes :

\medskip

\textbf{Sous-étape 1 :} Grâce au théorème de Banach-Steinhaus, on montre que
\[
\liminf_{j\to\infty} \liminf_{m\to\infty} p_{\Lambda(j)} (\mu_m ; \lambda_m) \geqslant p(\mu ; \lambda)
\]

\textbf{Sous-étape 2 :} On inverse ensuite les limites pour obtenir
\[
\liminf_{m\to\infty} p \big(\mu_{\Lambda(m+g(m)+\ell)}^{K_m} ; \lambda_m \big) \geqslant \limsup_{j\to\infty} \liminf_{m\to\infty} p_{\Lambda(j)} \big(\mu_{\Lambda(m+g(m)+\ell)}^{K_m} ; \lambda_m \big)
\]
Cela est vrai pour une certaine suite $(K_m)_{m \in \mathbb{N}}$ vérifiant bien la condition (\ref{compexh}). On utilise ici un argument d'uniformité dans la sous-additivité définissant les $p(\mu_{\Lambda(m+g(m)+\ell)} ; \lambda_m)$.

\medskip

\textbf{Sous-étape 1 :} Si $n \in \mathbb{N}$, le théorème de Banach-Steinhaus montre que les fonctions
\[
x \in K_n \mapsto \exp\langle \lambda_m | x \rangle
\]
convergent uniformément vers $\exp\langle \lambda | \cdot \rangle$ sur $K_n$. Ainsi, pour $m$ assez grand
\[
\int \exp \left\langle \lambda_m \bigg| \sum_{z \in \Lambda(j)} \xi(z) \right\rangle d \mu_m(\xi) \geqslant \int_{K_n^{\Lambda(j)}} \exp \left\langle \lambda \bigg| \sum_{z \in \Lambda(j)} \xi(z) \right\rangle d \mu_m(\xi) - \frac{1}{n}
\]
En passant à la limite inférieure en $m$, puis à la limite en $n$, on obtient\footnote{Pour $m \geqslant j \vee n$, on a
\[
\int_{\big\{ \xi_{\Lambda(j)} \in K_n^{\Lambda(j)} \big\}} \exp \left\langle \lambda \bigg| \sum_{z \in \Lambda(j)} \xi(z) \right\rangle d \mu_m(\xi) = \int_{\big\{ \xi_{\Lambda(j)} \in K_n^{\Lambda(j)} \big\}} \exp \left\langle \lambda \bigg| \sum_{z \in \Lambda(j)} \xi(z) \right\rangle d \mu(\xi)
\]}
\[
\liminf_{m\to\infty} \int \exp \left\langle \lambda_m \bigg| \sum_{z \in \Lambda(j)} \xi(z) \right\rangle d \mu_m(\xi) \geqslant \int \exp \left\langle \lambda \bigg| \sum_{z \in \Lambda(j)} \xi(z) \right\rangle d \mu(\xi)
\]
En prenant enfin la limite inférieure en $j$, on obtient
\[
\liminf_{j\to\infty} \liminf_{m\to\infty} p_{\Lambda(j)} (\mu_m ; \lambda_m) \geqslant p(\mu ; \lambda)
\]

\textbf{Sous-étape 2 :} Notons $\mu_m' = \mu_{\Lambda(m+g(m)+\ell)}$. On va se ramener à montrer que
\[
\liminf_{m\to\infty} p \big(\mu_m' ; \lambda_m \big) \geqslant \limsup_{j\to\infty} \liminf_{m\to\infty} p_{\Lambda(j)} \big(\mu_m' ; \lambda_m \big)
\]

En effet, on a d'une part, pour tout $j$,
\[
\liminf_{m\to\infty} p_{\Lambda(j)}(\mu_m' ; \lambda_m) \geqslant \liminf_{m\to\infty} p_{\Lambda(j)}(\mu_m ; \lambda_m)
\]
car
$$
\liminf_{m\to\infty} \frac{1}{|\Lambda(j)|} \log \mu_{\Lambda(m+g(m)+\ell)} \big( K_m^{\Lambda(m+g(m)+\ell)} \big) = 0
$$
D'autre part, en utilisant la proposition \ref{pmulambda},
\begin{align*}
p(\mu_m ; \lambda_m) = & \, \frac{1}{|\Lambda(m+g(m)+\ell)|} \log \int_{K_m^{\Lambda(m+g(m)+\ell)}} \exp \left\langle \lambda \bigg| \sum_{z \in \Lambda(m+g(m)+\ell)} \xi(z) \right\rangle d \mu(\xi)\\
 &\qquad - \frac{1}{|\Lambda(m+g(m)+\ell)|} \log \mu\big( K_m^{\Lambda(m+g(m)+\ell)} \big)
\end{align*}
et
\[
p(\mu_m' ; \lambda_m) = \frac{1}{|\Lambda(m+g(m)+\ell)|} \log \int \exp \left\langle \lambda \bigg| \sum_{z \in \Lambda(m+g(m)+\ell)} \xi(z) \right\rangle d \mu(\xi)
\]
Or, pour tout $m$,
\[
\lim_{n \to \infty} \int_{K_n^{\Lambda(m)}} \exp \left\langle \lambda \bigg| \sum_{z \in \Lambda(m)} \xi(z) \right\rangle d \mu(\xi) = \int \exp \left\langle \lambda \bigg| \sum_{z \in \Lambda(m)} \xi(z) \right\rangle d \mu(\xi)
\]
Ainsi, par un procédé diagonal, quitte à extraire une sous-suite de $(K_m)$, on obtient
\[
\liminf_{m \to \infty} p(\mu_m ; \lambda_m) = \liminf_{m \to \infty} p(\mu_m' ; \lambda_m)
\]
Reste donc à montrer que
\[
\liminf_{m\to\infty} p \big(\mu_m' ; \lambda_m \big) \geqslant \limsup_{j\to\infty} \liminf_{m\to\infty} p_{\Lambda(j)} \big(\mu_m' ; \lambda_m \big)
\]

Pour ce faire, on va montrer une forme d'uniformité sur l'inégalité sous-additive définissant les $p(\mu_m' ; \lambda_m)$ ; celle-ci est conséquence de la proposition suivante :

\begin{Pro}
Si $\mu$ est c.a.d.i. et c.l., et $m \in \mathbb{N}$, alors $\mu_{\Lambda(m)}$ est aussi c.a.d.i. et c.l., avec les mêmes paramètres.
\end{Pro}
\textbf{Démonstration :} Montrons par exemple le contrôle local. Notons $(t, \alpha)$ le paramètre de $\mu$ et fixons $z \in \mathbb{Z}^d$. On note $\Lambda_z$ l'unique boîte de la forme $a m + \Lambda(m)$ ($a \in \mathbb{Z}$) qui contient $z$. Alors, pour tout $\Lambda \subset \mathbb{Z}^d \setminus \{ z \}$,
\begin{align*}
\mu_{\Lambda(m)}&\, \big( \xi(z) \in t(V)\, V ; \xi_\Lambda \in C_\Lambda \big)\\
 &= \sum_i \mu_{\Lambda(m)} \big( \xi(z) \in t(V)\, V ; \xi_{\Lambda \cap \Lambda_z} \in C^i_{\Lambda \cap \Lambda_z} ; \xi_{\Lambda \setminus \Lambda_z} \in C^i_{\Lambda \setminus \Lambda_z} \big)\\
 &= \sum_i \mu_{\Lambda(m)} \big( \xi(z) \in t(V)\, V ; \xi_{\Lambda \cap \Lambda_z} \in C^i_{\Lambda \cap \Lambda_z}\big) \mu_{\Lambda(m)} \big(\xi_{\Lambda \setminus \Lambda_z} \in C^i_{\Lambda \setminus \Lambda_z} \big)\\
 &\geqslant \sum_i \alpha(V) \mu_{\Lambda(m)} \big( \xi_{\Lambda \cap \Lambda_z} \in C^i_{\Lambda \cap \Lambda_z}\big) \mu_{\Lambda(m)} \big(\xi_{\Lambda \setminus \Lambda_z} \in C^i_{\Lambda \setminus \Lambda_z} \big)\\
 &= \sum_i \alpha(V) \mu_{\Lambda(m)} \big( \xi_{\Lambda \cap \Lambda_z} \in C^i_{\Lambda \cap \Lambda_z} ; \xi_{\Lambda \setminus \Lambda_z} \in C^i_{\Lambda \setminus \Lambda_z} \big)\\
 &= \alpha(V) \mu_{\Lambda(m)}\, \big(\xi_\Lambda \in C_\Lambda \big)
\end{align*}
Le raisonnement est analogue pour le découplage.\qed

\medskip

Notons $\zeta_m$ un champ de loi $\mu_m'$. En vertu de la proposition précédente, ainsi que des propositions \ref{scadimos} et \ref{sclmos}, $\langle \lambda_m | \zeta_m \rangle$ est $(g, c)$-c.a.d.i. et $(t \circ \lambda_m^{-1}, \alpha \circ \lambda_m^{-1})$-c.l. Ainsi, l'inégalité (\ref{plsa}) se récrit
\begin{align*}
& \forall n \geqslant N(m)\nonumber\\
p_{\Lambda(n)}&(\mu_m' ; \lambda_m) \geqslant (1 - \rho_{j, n}) p_{\Lambda(j)}(\mu_m' ; \lambda_m)  - \frac{c(j)}{|\Lambda(j)|} - \rho_{j, n}\big(t(\lambda_m^{-1}(V_1)) - \log \alpha(\lambda_m^{-1}(V_1))\big)
\end{align*}

Plus précisément, comme $Y$ est tonnelé et comme les $\lambda_m$ sont bornées en tout point (étant donné que $\lambda_m \to \lambda$ faiblement), le tonneau
\[
T = \bigcap_{m=0}^{\infty} \lambda_m^{-1}\big( [-\textstyle{\frac{1}{2}}, {\frac{1}{2}}] \big)
\]
est un voisinage de $0$. Soit $V \in \mathcal{C}_0(Y)$ contenu dans $T$. La démonstration pour obtenir \ref{plsa} peut être refaite pour aboutir à
\begin{align*}
\forall n \geqslant &\, N(m)\nonumber\\
& p_{\Lambda(n)}(\mu_m' ; \lambda_m) \geqslant (1 - \rho_{j, n}) p_{\Lambda(j)}(\mu_m' ; \lambda_m) - \frac{c(j)}{|\Lambda(j)|} - \rho_{j, n}\big( t(V) - \log \alpha(V) \big)
\end{align*}

Prenant la limite en $n$, on obtient
\[
\forall m \in \mathbb{N} \quad \forall j \in \mathbb{N} \qquad p(\mu_m' ; \lambda_m) \geqslant p_{\Lambda(j)}(\mu_m' ; \lambda_m) - \frac{c'(j)}{|\Lambda(j)|}
\]
Puis, pour tout $\delta > 0$, il existe $j_\delta \in \mathbb{N}$ tel que
\[
\forall m \in \mathbb{N} \qquad p(\mu_m' ; \lambda_m) \geqslant \sup_{j\geqslant j_\delta} p_{\Lambda(j)}(\mu_m' ; \lambda_m) - \delta
\]
Alors
\begin{align*}
\liminf_{m\to\infty} p(\mu_m' ; \lambda_m) &\geqslant \liminf_{m\to\infty} \sup_{j\geqslant j_\delta} p_{\Lambda(j)}(\mu_m' ; \lambda_m) - \delta\\
&\geqslant \sup_{j\geqslant j_\delta} \liminf_{m\to\infty} p_{\Lambda(j)}(\mu_m' ; \lambda_m) - \delta\\
&\geqslant \limsup_{j\to\infty} \liminf_{m\to\infty} p_{\Lambda(j)}(\mu_m' ; \lambda_m) - \delta
\end{align*}
Prenant enfin la limite quand $\delta \to 0^+$, on obtient le résultat annoncé.\qed

\subsection{\'Egalité $s = -p^*$ : cas général}

Enfin, on peut énoncer :
\begin{Th}\label{spetoile}
Soient $\overleftarrow{\mathcal{X}} = (Y_i, f_i, f_{ij})_{i \leqslant j}$ un système projectif d'espaces de Fréchet séparables et $X$ sa limite projective. Soit $\sigma$ un champ à valeurs dans $\smash{\big(Y^{\mathbb{Z}^d}, \mathcal{B}^{\otimes \mathbb{Z}^d}\big)}$ tel que, pour tout $i \in J$, le champ $f_i(\sigma)$ soit c.a.d.i., c.l. et $R_\ell$-invariant. Alors,
\[
s(\sigma ; x) = - p^*(\sigma ; x)
\]
\end{Th}

\textbf{Démonstration :} Le théorème \ref{lspetoile} assure que
\[
\forall i \in I \quad \forall y_i \in Y_i \qquad s(f_i(\sigma) ; y_i) = - p^*(f_i(\sigma) ; y_i)
\]
Soit maintenant $x \in X$ ; on a :
\begin{align*}
s(\sigma ; x) &= \inf_{\substack{A \in \mathcal{C}\\ A \ni x}} \underline{s}(\sigma ; A)\\
 &= \inf_{i \in I} \left[ \inf_{\substack{V_i \in \mathcal{C}_i\\ V_i \ni f_i(x)}} \underline{s}(f_i(\sigma) ; V_i)\right]\\
 &= \inf_{i \in I} s\big(f_i(\sigma) ; f_i(x)\big)\\
 &= - \sup_{i \in I} p^*\big(f_i(\sigma) ; f_i(x)\big)
\end{align*}
la dernière égalité étant conséquence du théorème \ref{lspetoile}. Puis
\begin{align*}
\sup_{i \in I}\,\,& p^*\big(f_i(\sigma) ; f_i(x)\big)\\
 &= \sup_{i \in I} \sup_{\lambda_i \in Y_i^*} \left[ \langle \lambda_i \circ f_i | x \rangle - \lim_{n\to\infty} \frac{1}{|\Lambda(n)|} \log \int \exp\left\langle \lambda_i \circ f_i \bigg| \sum_{z \in \Lambda(n)} \sigma(z) \right\rangle d \mathbb{P} \right]\\
 &= \sup_{\lambda \in X^*} \left[ \langle \lambda | x \rangle - \lim_{n\to\infty} \frac{1}{|\Lambda(n)|} \log \int \exp\left\langle \lambda \bigg| \sum_{z \in \Lambda(n)} \sigma(z) \right\rangle d \mathbb{P} \right]\\
 &= p^*(\sigma ; x)
\end{align*}
D'où le résultat attendu.\qed

\bibliographystyle{alpha-fr}
\bibliography{../cramer}
%\bibliography{cramer}

\end{document}